\documentclass[10pt]{amsart}
\usepackage{mathptmx}
\usepackage{amsmath}
\usepackage{amssymb}
\usepackage{array}
\usepackage{geometry}
\usepackage[bookmarks=true,colorlinks=true, pdfstartview=FitV, linkcolor=black, citecolor=blue, urlcolor=black]{hyperref}

\usepackage{color}
\definecolor{DarkRed}{rgb}{0.55,.00,0.2}
\definecolor{DarkGrey}{rgb}{0.35,.35,0.35}

\theoremstyle{definition}

\theoremstyle{remark}

\numberwithin{equation}{section}

\begin{document}

\title{ On the Weber integral equation \\ and solution to the Weber-Titchmarsh problem}

\author{S. Yakubovich}

\address{Department of Mathematics, Fac. Sciences of University of Porto,Rua do Campo Alegre,  687; 4169-007 Porto (Portugal)}

\email{ syakubov@fc.up.pt}

\begin{abstract}  We derive sufficient conditions for the existence of the Weber formal solution of the corresponding integral equation, related to the familiar Weber-Orr integral transforms.  This gives a solution to the old Weber-Titchmarsh problem (posed in {\it Proc. Lond. Math. Soc.}  {\bf 22}(2) (1924),  pp.15, 16.)   Our method involves properties of the inverse Mellin transform  of integrable functions.  The Mellin-Parseval equality and some integrals with the associated Legendre functions are used. 
\end{abstract}

\keywords{Keywords  here}

\subjclass[2000]{Primary 44A15, 44A35,  33C10;    Secondary 33C05, 33C45 }

\date{\today}

\keywords{  Weber-Orr integral transforms,   Mellin transform, Bessel functions, Associated Legendre functions}

\maketitle

\markright{\rm \centerline{SOLUTION TO THE WEBER-TITCHMARSH PROBLEM}}

\section{Introduction and preliminary results}

In \cite{titweb} E.C. Titchmarsh formally showed  that an arbitrary complex-valued function $g(x),\ x \in \mathbb{R}_+$ can be expanded in terms of  the following repeated  integral 

$$g(x)=  {x\over  J_\nu^2(ax)+ Y_\nu^2(ax)} \int_a^\infty C_\nu (xt, xa) t \int_0^\infty   C_\nu (t\xi, a\xi) g(\xi) d\xi dt,\eqno(1)$$  
where $a >0$,  $\nu \in \mathbb{C}$, $J_\nu(z), Y_\nu(z)$ are Bessel functions of the first and second kind \cite{erd}, Vol. II and 

$$C_\nu(\alpha, \beta)=  J_\nu(\alpha) Y_\nu(\beta) -   Y_\nu(\alpha) J_\nu(\beta).\eqno(2)$$
Expansion (1) is related to the familiar Weber-Orr integral expansions of an arbitrary function $f(x)$ as repeated integrals

$$f(x)=  \int_0^\infty  {t \  C_\nu (xt, at)\over  J_\nu^2(at)+ Y_\nu^2(at)}\int_a^\infty   C_\nu (\xi t, at) \xi f(\xi) d\xi dt,\eqno(3)$$  

$$f(x)=   \int_a^\infty C_\nu (xt, xa) t \int_0^\infty    { C_\nu (t \xi, a\xi) \over  J_\nu^2(a\xi)+ Y_\nu^2(a\xi)} \xi f (\xi) d\xi dt,\eqno(4)$$  
which are different from (1).   Combining with (3), (4), Titchmarsh proved formally (1) for $a=1$ (see \cite{titweb}, p. 15).  He posed the problem to find sufficient conditions for the validity of expansion (1) in order to  solve the following Weber integral equation with respect to $g$

$$f(x)= \int_0^\infty   C_\nu (x \xi, a\xi) g(\xi) d\xi,\eqno(5)$$
where $f(x),\  x  \in \mathbb{R}_+$  is a given function.  As far as  the author is aware this question is still open.    Our method will be based on the use of the Mellin transform \cite{tit}.  Precisely, the Mellin transform is defined in  $L_{\mu, p}(\mathbb{R}_+),\ 1 < p \le 2$  by the integral  
$$f^*(s)= \int_0^\infty f(x) x^{s-1} dx,\eqno(6)$$
 being convergent  in mean with respect to the norm in $L_q(\mu- i\infty, \mu + i\infty),\   q=p/(p-1)$.   Moreover, the  Parseval equality holds for $f \in L_{\mu, p}(\mathbb{R}_+),\  g \in L_{1-\mu, q}(\mathbb{R}_+)$
$$\int_0^\infty f(x) g(x) dx= {1\over 2\pi i} \int_{\mu- i\infty}^{\mu+i\infty} f^*(s) g^*(1-s) ds.\eqno(7)$$
The inverse Mellin transform is given accordingly
 $$f(x)= {1\over 2\pi i}  \int_{\mu- i\infty}^{\mu+i\infty} f^*(s)  x^{-s} ds,\eqno(8)$$
where the integral converges in mean with respect to the norm  in   $L_{\mu, p}(\mathbb{R}_+)$
$$||f||_{\mu,p} = \left( \int_0^\infty  |f(x)|^p x^{\mu p-1} dx\right)^{1/p}.\eqno(9)$$
In particular, letting $\mu= 1/p$ we get the usual space $L_1(\mathbb{R}_+)$.   A special class of functions related to the Mellin transform (6) and its inversion (8),  was introduced in \cite{class}. Indeed, we have

{\bf Definition 1 (\cite{class})}. Denote by ${\mathcal M}^{-1}(L_c)$ the space of
functions $f(x), \ x \in \mathbb{R}_+$,  representable by inverse
Mellin transform (8) of integrable functions $f^{*}(s) \in L_{1}(c)$ on
the vertical line $c =\{s \in \mathbb{C}:  \mu ={\rm Re s} = c_0\}$.

The  space ${\mathcal M}^{-1}(L_c)$  with  the  usual operations  of
addition   and multiplication by scalar is a linear vector space. If
the norm in ${\mathcal M}^{-1}(L_c)$ is introduced by the formula
$$ \big\vert\big\vert f \big\vert\big\vert_{{\mathcal
M}^{-1}(L_c)}= {1\over 2\pi }\int^{+\infty}_{-\infty} |
f^{*}\left(c_0 +it\right)| dt,\eqno(10)$$
then it becomes  a Banach space.

 {\bf Definition 2 (\cite{class}, \cite{yal}}).  Let $\mu\neq 0,\ c_1, c_2 \in \mathbb{R}$ be such that $2 \hbox{sign}\ c_1 + \hbox{sign}\  c_2 \ge 0$. By ${\mathcal M}_{c_1,c_2}^{-1}(L_c)$ we denote the space of functions $f(x), x \in \mathbb{R}_+$, representable in the form (8), where $s^{c_2}e^{\pi c_1|s|} f^*(s) \in L_1(c)$.

It is a Banach space with the norm
$$ \big\vert\big\vert f \big\vert\big\vert_{{\mathcal
M}_{c_1,c_2}^{-1}(L_c)}= {1\over 2\pi }\int_{c} e^{\pi c_1|s|}
|s^{c_2} f^{*}(s) ds|.$$
In particular, letting $c_1=c_2=0$ we get the space ${\mathcal M}^{-1}(L_c)$. Moreover, it is easily seen the inclusion

$${\mathcal M}_{d_1,d_2}^{-1}(L_c) \subseteq {\mathcal
M}_{c_1,c_2}^{-1}(L_c)$$ when $2 \hbox{sign}(d_1- c_1) + \hbox{sign}
(d_2-c_2) \ge 0$.

\section{Solution to the Weber-Titchmarsh problem}

The goal  of this Note is to prove the following 

{\bf Theorem.} {\it Let $a > 0, \ \nu \in \mathbb{C}, \   0< {\rm Re } \nu < 1/2,  g(x)  \in  {\mathcal M}_{0,1} ^{-1}(L_c)$ with $c =\{s \in \mathbb{C}:  -1< {\rm Re} s < 0\}$.   Then for almost all $x >0$ expansion $(1)$ holds, where the inner and outer integrals are understood in the  improper sense.} 

\begin{proof}    We begin,   writing $g(\xi)$ in terms of the inverse Mellin transform (8) of the reciprocal function $g^*(s) \in L_1(c)$. Then substituting this expression into (5), we change  the order of integration when $\xi \in [0, N]$ to find

$$f(x)=  {1\over 2\pi i}  \lim_{N\to \infty} \int_{\mu- i\infty}^{\mu+i\infty} g^*(s) \int_0^N   C_\nu (x \xi, a\xi)  \xi^{-s}  d\xi ds,\eqno(11) $$
where the interchange is ensured by the Fubini theorem via absolute convergence of the repeated integral. Indeed, it  is due to conditions of the theorem and asymptotic behavior of Bessel functions at infinity and near the origin \cite{erd} via the assumption  $-1<  \mu < 0.$   Then kernel (2) behaves as follows

$$ C_\nu (x \xi, a\xi)  = O( 1),\ \xi \to 0+,\ \nu \neq  0,$$ 

$$ C_\nu (x \xi, a\xi)  = O( \log \xi),\ \xi \to 0+,\  \nu = 0,$$ 

$$ C_\nu (x \xi, a\xi)  = - {2\over \pi \xi \sqrt {a x}} \left[ \sin(\xi (x-a)) + O\left( {1\over \xi} \right)\right],\  \xi \to \infty.$$  
Hence, for fixed $x, a, N$ 

$$\int_{\mu- i\infty}^{\mu+i\infty} \left| g^*(s) \right|\int_0^N  \left|  C_\nu (x \xi, a\xi)\right|  \xi^{-\mu}  d\xi |ds| < \infty.$$
In the meantime, the integral with respect to $\xi$ over $\mathbb{R}_+$ can be calculated with the use of relation (2.13.15.4) in \cite{prud}, Vol. 2,   Boltz's and self-transformation formulae  for the Gauss hypergeometric function and basic relations between the associated Legendre functions (see details in \cite{erd}, Vol. I).  Thus we obtain

$$\int_0^\infty   C_\nu (x \xi, a\xi)  \xi^{-s}  d\xi  =   {2^{1 -s}\over \pi} \ e^{i\nu\pi} \left(x^2- a^2\right)^{(s-1)/2} 
\frac{\Gamma((1-s)/2)}{ \Gamma( (1- 2\nu + s)/2)} Q^{-\nu}_{(s-1)/2} \left( {x^2+a^2 \over x^2-a^2}\right),\eqno(12)$$
where $x  > a > 0, \Gamma(z)$ is Euler's gamma-function and $Q^\nu_\mu (z)$ is the associated Legendre function of the second kind \cite{erd}, Vol. I, \cite{vira}.   Hence, we write (11) in the form 

$$ f(x)=  { e^{i\nu\pi} \over \pi^2  i}  \int_{\mu- i\infty}^{\mu+i\infty} g^*(s) \frac{\Gamma((1-s)/2)}{ \Gamma( (1- 2\nu + s)/2)} Q^{-\nu}_{(s-1)/2} \left( {x^2+a^2 \over x^2-a^2}\right) 2^{-s}  \left(x^2- a^2\right)^{(s-1)/2} ds$$

$$-  {1\over 2\pi i}  \lim_{N\to \infty} \int_{\mu- i\infty}^{\mu+i\infty} g^*(s) \int_N^\infty   C_\nu (x \xi, a\xi)  \xi^{-s}  d\xi ds\eqno(13)$$
and will prove that 

$$\lim_{N\to \infty} \int_{\mu- i\infty}^{\mu+i\infty} g^*(s) \int_N^\infty   C_\nu (x \xi, a\xi)  \xi^{-s}  d\xi ds= 0,\ -1< \mu < 0.$$
To do this, we will need the exact asymptotic behavior at infinity of the kernel (2) (see above).  Then, substituting the main asymptotic term  and integrating by parts in the integral by $\xi$, we deduce for fixed $x > a$

$$  - {2\over \pi \sqrt{a x} } \int_N^\infty    \sin(\xi (x-a))  \xi^{-s-1}  d\xi =  -   {2 \cos(N (x-a)) \over \pi (x-a) \sqrt{a x} }   N^{-s-1} +  O( (s+1) N^{-s-1} ), \ N \to \infty.$$
Hence, under conditions of the theorem 

$$ \left| \int_{\mu- i\infty}^{\mu+i\infty} g^*(s) \int_N^\infty   C_\nu (x \xi, a\xi)  \xi^{-s}  d\xi ds\right| \le C\  N^{-\mu-1}  \int_{\mu- i\infty}^{\mu+i\infty} \left| g^*(s) \right| (|s| +2)|ds| \to 0,\ N \to \infty,  $$
where $C >0$ is an absolute constant.   We note that will use the same notation below for different positive constants. Now, taking into account (13), let us consider the following sequence of functions (see (1)) 

 $$G_N(x)=     { e^{i\nu\pi} \over \pi^2  i}  \int_a^N  C_\nu (xt, xa) t  \int_{\mu- i\infty}^{\mu+i\infty} g^*(s) 
\frac{\Gamma((1-s)/2)}{ \Gamma( (1- 2\nu + s)/2)} $$

$$\times Q^{-\nu}_{(s-1)/2} \left( {t^2+a^2 \over t^2-a^2}\right) 2^{-s}  \left(t^2- a^2\right)^{(s-1)/2} ds dt.\eqno(14)$$
Choosing a complex variable $w$ from the vertical strip $-1 < {\rm Re} w < - 1/2$, we multiply both sides of (14) by $x^w$ and  integrate with respect to $x$ over $\mathbb{R}_+$. Hence, changing  the order of integration in its right-hand side by Fubini's theorem owing to  the same motivation as above when $N$ is fixed.  Hence,  it becomes

$$\int_0^\infty G_N(x) x^w dx =      { e^{i\nu\pi} \over \pi^2  i}  \int_a^N   t \int_0^\infty   C_\nu (xt, xa) x^w dx  \int_{\mu- i\infty}^{\mu+i\infty} g^*(s)  \frac{\Gamma((1-s)/2)}{ \Gamma( (1- 2\nu + s)/2)} $$

$$\times Q^{-\nu}_{(s-1)/2} \left( {t^2+a^2 \over t^2-a^2}\right) 2^{-s}  \left(t^2- a^2\right)^{(s-1)/2} ds dt.$$
Making use (12), we find the equality

$$\int_0^\infty G_N(x) x^w dx =     { 2^{1+w} e^{2 i\nu\pi} \over \pi^3  i}  \frac{\Gamma((1+w)/2)}{ \Gamma( (1- 2\nu -w)/2)}\int_{\mu- i\infty}^{\mu+i\infty} g^*(s)   \frac{\Gamma((1-s)/2)}{ \Gamma( (1- 2\nu + s)/2)} 2^{-s} $$

$$\times   \int_a^N   t   \left(t^2- a^2\right)^{(s- w-2)/2} Q^{-\nu}_{- (w+1)/2} \left( {t^2+a^2 \over t^2-a^2}\right) Q^{-\nu}_{(s-1)/2} \left( {t^2+a^2 \over t^2-a^2}\right)   dt ds.\eqno(16)$$
Then  due to elementary substitutions the inner integral with respect to $t$ in (16) can be written as
$$\int_a^N   t   \left(t^2- a^2\right)^{(s- w-2)/2} Q^{-\nu}_{- (w+1)/2} \left( {t^2+a^2 \over t^2-a^2}\right) Q^{-\nu}_{(s-1)/2} \left( {t^2+a^2 \over t^2-a^2}\right)   dt $$
$$= {1\over 2} a^{s-w}  \int_1^\infty   \left(t- 1\right)^{(s- w-2)/2} Q^{-\nu}_{- (w+1)/2} \left( {t+1 \over t-1}\right) Q^{-\nu}_{(s-1)/2} \left( {t+1 \over t- 1}\right)   dt  -  I_N(s,w), $$
where

$$ I_N(s,w) =  {1\over 2}  a^{s-w}  \int_{(N/a)^2}^\infty    \left(t- 1\right)^{(s- w-2)/2} Q^{-\nu}_{- (w+1)/2} \left( {t+1 \over t-1}\right) Q^{-\nu}_{(s-1)/2} \left( {t+1 \over t-1}\right)   dt.\eqno(17) $$
Meanwhile,  Parseval equality (7) and relation (8.4.42.5) in \cite{prud}, Vol. 3 permit to get the equality

$${1\over 2} a^{s-w}  \int_1^\infty   \left(t- 1\right)^{(s- w-2)/2} Q^{-\nu}_{- (w+1)/2} \left( {t+1 \over t-1}\right) Q^{-\nu}_{(s-1)/2} \left( {t+1 \over t- 1}\right)   dt $$

$$=  {e^{-2i\nu \pi} \over 16\pi i}\  a^{s-w}  \Gamma( (1+ s)/2) \Gamma( (1- w)/2)  \Gamma( (1- 2\nu+ s)/2) \Gamma( (1- 2\nu - w)/2)$$

$$\times \int_{\gamma - i\infty}^{\gamma +i\infty} \frac {\Gamma( (\nu +w +1)/2 - \tau) \Gamma( (w+1-\nu)/2 -\tau ) \Gamma( (\nu -s -1)/2 + \tau) \Gamma(\tau-   (1+ s + \nu )/2 )} {\Gamma( 1+ \nu/2 -\tau ) \Gamma( 1- \nu/2 -\tau)\Gamma( \tau +\nu/2 ) \Gamma( \tau - \nu/2 )} d\tau,$$
where $-1 <   {\rm Re} s < 0, \  -1  <  {\rm Re} w < -1/2$  and 
$$ {1\over 2} \left( 1+ {\rm Re} (s + \nu) \right) < \gamma <    {1\over 2} \left( 1+  {\rm Re} (w  - \nu) \right).$$
However, the latter integral can be calculated as the sum of residues at the left-hand poles of gamma-functions with the use of Slater's theorem \cite{mar},  and  values of the hypergeometric function ${}_2F_1$ at the unity \cite{erd}, Vol. I, \cite{prud}, Vol. 3.  Thus   we obtain the result

$${1\over 2\pi i} \int_{\gamma - i\infty}^{\gamma +i\infty} \frac {\Gamma( (\nu +w +1)/2 - \tau) \Gamma( (w+1-\nu)/2 -\tau ) \Gamma( (\nu -s -1)/2 + \tau) \Gamma(\tau-   (1+ s + \nu )/2 )} {\Gamma( 1+ \nu/2 -\tau ) \Gamma( 1- \nu/2 -\tau)\Gamma( \tau +\nu/2 ) \Gamma( \tau - \nu/2 )} d\tau$$

$$= 2^{s-w+1} {\cos(\pi\nu)\over \sqrt\pi } \frac {\Gamma( (1+s- w)/2) \Gamma( (w-s)/2) \Gamma( (w -s)/2 -\nu) \Gamma((w-s)/2+ \nu )} {\Gamma( (1+ s)/2 ) \Gamma( (1-s)/2 )\Gamma((1+ w)/2 ) \Gamma( (1-  w)/2)},$$
which leads after simple substitution to the value of possibly new integral with the product of the associated Legendre functions of the second kind

$$\int_1^\infty   \left(t- 1\right)^{(w- s)/2 - 1} Q^{-\nu}_{- (w+1)/2} \left( t\right) Q^{-\nu}_{(s-1)/2} \left( t\right)   dt =  2^{(s-w)/2 -1}  e^{-2i\nu \pi}  {\cos(\pi\nu)\over \sqrt\pi } $$

$$\times \frac {\Gamma( (1+s- w)/2) \Gamma( (w-s)/2) \Gamma( (w -s)/2 -\nu) \Gamma((w-s)/2+ \nu )} { \Gamma( (1-s)/2 )\Gamma((1+ w)/2 ) }$$

$$\times  \Gamma( (1+ s)/2 - \nu) \Gamma( (1-w)/2- \nu),   {\rm Re}  \nu <  {\rm Re} \left( {w -s\over 2} \right) < {1\over 2},\eqno(18)$$
when  $w$ is related to $s$ by the condition $- 1/2 >  {\rm Re} w >  \hbox{max} \left( 2{\rm Re} \nu + {\rm Re} s, -1\right) $.
Further, recalling (16) and substituting the value of the integral (18), we write it in the form

$$\int_0^\infty G_N(x) x^w dx =     { \cos(\pi\nu)  \over 2 \pi^3 \sqrt \pi  i} \int_{\mu- i\infty}^{\mu+i\infty} g^*(s)   \Gamma( (1+s- w)/2) \Gamma( (w-s)/2)$$

$$\times  \Gamma( (w -s)/2 -\nu) \Gamma((w-s)/2+ \nu ) a^{s-w} ds $$

$$-   { 2^{1+w} e^{2 i\nu\pi} \over \pi^3  i}  \frac{\Gamma((1+w)/2)}{ \Gamma( (1- 2\nu -w)/2)}\int_{\mu- i\infty}^{\mu+i\infty} g^*(s)   \frac{\Gamma((1-s)/2)}{ \Gamma( (1- 2\nu + s)/2)}  I_N(s,w)  2^{-s} ds,\eqno(19) $$
where $I_N(s,w)$ is defined by (17).  Meanwhile, the first term in the right-hand side of (19) can be treated by the Parseval equality (7) for the Mellin transform with the use of relation (8.4.20.35) in \cite{prud}, Vol. 3, which gives the Mellin- Barnes integral representation for the kernel  $J_\nu^2(ax)+ Y_\nu^2(ax)$, namely

$$J_\nu^2(ax)+ Y_\nu^2(ax) =   { \cos(\pi\nu)  \over 2 \pi^3 \sqrt \pi  i} \int_{\gamma - i\infty}^{\gamma +i\infty}  \Gamma( (1-s)/2) \Gamma( s/2)  \Gamma( s/2 -\nu) \Gamma(s/2+ \nu ) (ax)^{-s} ds,\   2  {\rm Re}  \nu < \gamma <  1. $$
Therefore, under conditions of the theorem we obtain 

$$  { \cos(\pi\nu)  \over 2 \pi^3 \sqrt \pi  i} \int_{\mu- i\infty}^{\mu+i\infty} g^*(s)   \Gamma( (1+s- w)/2) \Gamma( (w-s)/2)$$

$$\times  \Gamma( (w -s)/2 -\nu) \Gamma((w-s)/2+ \nu ) a^{s-w} ds  =  \int_0^\infty g(x) \left[ J_\nu^2(ax)+ Y_\nu^2(ax) \right] x^{w-1} dx.\eqno(20)$$
On the other hand,    making use the same substitution, we write (17) as follows 

$$ I_N(s,w) =   2^{s-w-1}  a^{s-w}  \int_{1}^{(N^2+a^2)/ (N^2- a^2)}     \left(t- 1\right)^{(w-s)/2- 1} Q^{-\nu}_{- (w+1)/2} \left( t\right) Q^{-\nu}_{(s-1)/2} \left( t\right)   dt.\eqno(21)$$
The associated Legendre functions of the second kind  $Q^{-\nu}_{- (w+1)/2} \left( t\right),\  Q^{-\nu}_{(s-1)/2} \left( t\right)$ , in turn, can be represented in terms of the Euler integral for the Gauss hypergeometric function \cite{erd}, Vol. I and we find

$$Q^{-\nu}_{- (w+1)/2} \left( t\right) =    \frac{ 2^{(w-1)/2}   e^{-i\nu\pi} \sqrt\pi\  \Gamma((1-w - 2\nu)/2)}{ \Gamma( (3- w - 2\nu)/4) \Gamma((1-w+2\nu)/4)}\  (t^2-1)^{-\nu/2} $$

$$\times \int_0^1u^{- (w+1+2\nu)/4} (1-u)^{ (2\nu- w-3)/4} (t^2-u)^{ (2\nu+ w-  1)/4} du,\eqno(22)$$

$$Q^{-\nu}_{ (s-1)/2} \left( t\right) =    \frac{ 2^{-(s+1)/2}   e^{-i\nu\pi} \sqrt\pi\  \Gamma((s+1 - 2\nu)/2)}{ \Gamma((3+ s  -2 \nu)/4) \Gamma((1+s+2\nu)/4)}\  (t^2-1)^{-\nu/2} $$

$$\times \int_0^1u^{- (1-s +2\nu)/4} (1-u)^{ (s-3+2\nu)/4} (t^2-u)^{ (2\nu-  1-s)/4} du, \  {\rm Re} s >  2 {\rm Re}  \nu-1.\eqno(23)$$
Hence,  we obtain the following uniform estimates for the associated Legendre functions of the second kind 

$$\left| Q^{-\nu}_{- (w+1)/2} \left( t\right) \right| \le \sqrt\pi  \ 2^{( {\rm Re} w-1)/2}  (t^2-1)^{- {\rm Re} \nu/2}  \left| \frac{ \Gamma((1-w - 2\nu)/2)}{ \Gamma( (3- w - 2\nu)/4) \Gamma((1-w+2\nu)/4)}\right|$$

$$\times \int_0^1u^{- {\rm Re} (w+1+2\nu)/4} (1-u)^{ {\rm Re}\nu -1} du = O \left( |w|^{ - {\rm Re} \nu }  (t^2-1)^{- {\rm Re} \nu/2} \right),\ t > 1,\eqno(24)$$

$$\left| Q^{-\nu}_{ (s-1)/2} \left( t\right)\right| \le \sqrt\pi\  2^{-({\rm Re} s+1)/2}  (t^2-1)^{-{\rm Re} \nu/2} \left| \frac{  \Gamma((s+1 - 2\nu)/2)}{ \Gamma((3+ s  -2 \nu)/4) \Gamma((1+s+2\nu)/4)}\right|\   $$

$$\times \int_0^1u^{-  {\rm Re} (1-s +2\nu)/4} (1-u)^{ {\rm Re}\nu -1} du=  O \left( |s|^{ - {\rm Re} \nu }   (t^2-1)^{- {\rm Re} \nu/2} \right),\ t >1,\eqno(25)$$
since via Stirling's asymptotic formula for the gamma-function \cite{erd}, Vol. I

$$ \left| \frac{ \Gamma((1-w - 2\nu)/2)}{ \Gamma( (3- w - 2\nu)/4) \Gamma((1-w+2\nu)/4)}\right| =  O \left( |w|^{ - {\rm Re} \nu }\right),\  \left|{\rm Im} w\right| \to \infty,$$ 

$$ \left| \frac{  \Gamma((s+1 - 2\nu)/2)}{ \Gamma((3+ s  -2 \nu)/4) \Gamma((1+s+2\nu)/4)}\right| =  O \left( |s|^{ - {\rm Re} \nu }\right),\  \left|{\rm Im} s\right| \to \infty.$$ 
Therefore, returning to (21), we have by the straightforward estimate

$$\left| I_N(s,w) \right| \le C  |w s|^{ - {\rm Re} \nu }   \int_{1}^{(N^2+a^2)/ (N^2- a^2)}     \left(t- 1\right)^{{\rm Re} (w-s-2\nu)/2- 1} dt$$

$$  =  O\left(   |w s|^{ - {\rm Re} \nu }   N^{ {\rm Re} (s+2\nu- w)}\right),\   N\to \infty\eqno(26)$$
under condition (see (18))  ${\rm Re} w >  2{\rm Re} \nu + {\rm Re} s$.  

Further, substituting (20) into (19), using the obtained estimates and Stirling's asymptotic formula for the gamma-function, we find

$$\left| \int_0^\infty \left[ G_N(x)-  g(x) \left[ J_\nu^2(ax)+ Y_\nu^2(ax)\right] x^{-1}\right]  x^w dx \right| $$

$$\le   C    |w |^{ - {\rm Re} \nu }   N^{ {\rm Re} (s+2\nu- w)} \left| \frac{\Gamma((1+w)/2)}{ \Gamma( (1- 2\nu -w)/2)}\right| \int_{\mu- i\infty}^{\mu+i\infty} \left| g^*(s) \right|  |s|^{ - {\rm Re} \nu }  \left|\frac{\Gamma((1-s)/2)}{ \Gamma( (1- 2\nu + s)/2)}  ds\right|$$

$$\le C |w |^{  {\rm Re} w} N^{ {\rm Re} (s+2\nu- w)}  \int_{\mu- i\infty}^{\mu+i\infty} \left| g^*(s) \right|  |s|^{ - \mu }  \left| ds\right| = O\left( N^{ {\rm Re} (s+2\nu- w)}  \right),\ N \to \infty.$$
This is because $g \in   {\mathcal M}_{0,1} ^{-1}(L_c)$.   Thus  for all $w$ from the strip $-1/2  >  {\rm Re} w > \mu + 2 {\rm Re} \nu$ it yields the equality 

$$\lim_{N\to \infty}   \int_0^\infty  G_N(x)  x^w dx =  \int_0^\infty  g(x) \left[ J_\nu^2(ax)+ Y_\nu^2(ax)\right] x^{w -1} dx.\eqno(27)$$
The latter equality means that the sequence of Mellin transforms (6) of variable $w$ of functions $x G_N(x)$ converges pointwisely to the Mellin transform of $g(x) \left[ J_\nu^2(ax)+ Y_\nu^2(ax)\right]$.  On the other hand, recalling (19), (20) and the generalized Minkowski inequality, we estimate $L_2$-norm of the difference of these functions. In fact, using (26), we have

$$\left( \int_{-\infty}^\infty \left|  \int_0^\infty  \left[ x G_N(x) -  g(x) \left[ J_\nu^2(ax)+ Y_\nu^2(ax)\right] \right] x^{{\rm Re w} +i\tau-1} dx\right|^2 d\tau \right)^{1/2} $$

$$= {2^{1+{\rm Re} w}\over \pi^{3} } \left( \int_{-\infty}^\infty \left|  \frac{\Gamma((1+ {\rm Re} w +i\tau )/2)}{ \Gamma( (1- 2\nu - {\rm Re} w - i\tau)/2)}\int_{\mu- i\infty}^{\mu+i\infty} g^*(s)   \frac{\Gamma((1-s)/2)}{ \Gamma( (1- 2\nu + s)/2)}  I_N(s,w)  2^{-s} ds \right|^2 d\tau \right)^{1/2} $$

$$\le {2^{1+{\rm Re} w- \mu}\over \pi^{3} } \int_{\mu- i\infty}^{\mu+i\infty}  \left| g^*(s)   \frac{\Gamma((1-s)/2)}{ \Gamma( (1- 2\nu + s)/2)}\right| \left( \int_{-\infty}^\infty  \left| \frac{\Gamma((1+ {\rm Re} w +i\tau )/2)}{ \Gamma( (1- 2\nu - {\rm Re} w - i\tau)/2)}\right|^2 \ \left| I_N(s,w) \right|^2 d\tau \right)^{1/2} | ds|$$

$$\le C \  N^{ \mu+ {\rm Re} (2\nu- w)}  \int_{\mu- i\infty}^{\mu+i\infty}  \left| g^*(s) \right|  |s|^{-\mu} |ds|  \left( \int_{-\infty}^\infty  \left( ({\rm Re} w)^2  + \tau^2 \right)^{{\rm Re} w}  d\tau \right)^{1/2} =  O\left(   N^{ \mu+ {\rm Re} (2\nu- w)}\right) \to 0, \  N \to \infty$$
under conditions of the theorem and the choice of $w$ from the strip  $-1/2  >  {\rm Re} w > \mu + 2 {\rm Re} \nu$.  This estimates imply that the limit (27) exists also in the mean square sense.  However, an analog of the Plancherel theorem for the Mellin transform (see \cite{tit}, Th. 71 ) and Parseval equality (7) say 

$${1\over 2\pi} \int_{-\infty}^\infty \left|  \int_0^\infty  \left[ x G_N(x) -  g(x) \left[ J_\nu^2(ax)+ Y_\nu^2(ax)\right] \right] x^{{\rm Re w} +i\tau-1} dx\right|^2 d\tau$$

$$ =  \int_{0}^\infty \left|  x G_N(x) -  g(x) \left[ J_\nu^2(ax)+ Y_\nu^2(ax)\right] \right|^2  x^{ 2 {\rm Re w} -1} dx. $$ 
Consequently, the right-hand side of the latter equality tends to zero as well when $N \to \infty$ and we get the value of the limit 

$$\lim_{N\to \infty}  x  G_N(x)  =    g(x) \left[ J_\nu^2(ax)+ Y_\nu^2(ax)\right],\ x >0, \eqno(28)$$
in the mean square sense.  Moreover, as it is known,  this limit holds for almost all $x >0$ for some subsequence $G_{N_k}$, namely,

$$\lim_{k\to \infty}   G_{N_k} (x)  =    g(x)\   x^{-1} \left[ J_\nu^2(ax)+ Y_\nu^2(ax)\right].$$
 Therefore, in order to complete the proof, we need to show that the sequence $G_N$ is a Cauchy one.  Indeed, taking (14), we write for some positive big enough $M, N, \  M > N$ and fixed $x >0$ 

$$G_M(x)- G_N(x) =     { e^{i\nu\pi} \over \pi^2  i}  \int_N^M   C_\nu (xt, xa) t  \int_{\mu- i\infty}^{\mu+i\infty} g^*(s) 
\frac{\Gamma((1-s)/2)}{ \Gamma( (1- 2\nu + s)/2)} $$

$$\times Q^{-\nu}_{(s-1)/2} \left( {t^2+a^2 \over t^2-a^2}\right) 2^{-s}  \left(t^2- a^2\right)^{(s-1)/2} ds dt.$$
Hence, appealing to (2), the asymptotic behavior of Bessel functions of the first and second kind at infinity and the estimate (25), we obtain

$$G_M(x)- G_N(x) =     { e^{i\nu\pi} \sqrt 2\over \pi^2 \sqrt {\pi a}  i} \left[  Y_\nu(xa)  \int_N^M  \cos (xt - \pi(2\nu+1)/4) \left(\sqrt t + O\left(t^{-1/2} \right)\right)\right. $$

$$\times \int_{\mu- i\infty}^{\mu+i\infty} g^*(s) \frac{\Gamma((1-s)/2)}{ \Gamma( (1- 2\nu + s)/2)} Q^{-\nu}_{(s-1)/2} \left( {t^2+a^2 \over t^2-a^2}\right) 2^{-s}  \left(t^2- a^2\right)^{(s-1)/2} ds dt$$

$$- J_\nu(xa)  \int_N^M  \sin  (xt - \pi(2\nu+1)/4) \left(\sqrt t + O\left(t^{-1/2} \right)\right) $$

$$\left. \times \int_{\mu- i\infty}^{\mu+i\infty} g^*(s) \frac{\Gamma((1-s)/2)}{ \Gamma( (1- 2\nu + s)/2)} Q^{-\nu}_{(s-1)/2} \left( {t^2+a^2 \over t^2-a^2}\right) 2^{-s}  \left(t^2- a^2\right)^{(s-1)/2} ds dt\right] $$

$$=  { e^{i\nu\pi} \sqrt 2\over \pi^2 \sqrt {\pi a}  i} \left[  Y_\nu(xa)  \int_N^M  \cos (xt - \pi(2\nu+1)/4) \sqrt t \right. $$

$$\times \int_{\mu- i\infty}^{\mu+i\infty} g^*(s) \frac{\Gamma((1-s)/2)}{ \Gamma( (1- 2\nu + s)/2)} Q^{-\nu}_{(s-1)/2} \left( {t^2+a^2 \over t^2-a^2}\right) 2^{-s}  \left(t^2- a^2\right)^{(s-1)/2} ds dt$$

$$- J_\nu(xa)  \int_N^M  \sin  (xt - \pi(2\nu+1)/4) \sqrt t   \int_{\mu- i\infty}^{\mu+i\infty} g^*(s) \frac{\Gamma((1-s)/2)}{ \Gamma( (1- 2\nu + s)/2)} $$

$$\left. \times Q^{-\nu}_{(s-1)/2} \left( {t^2+a^2 \over t^2-a^2}\right) 2^{-s}  \left(t^2- a^2\right)^{(s-1)/2} ds dt\right]+  O\left( N^{\mu+ {\rm Re}\nu- 1/2}\right). \eqno(29)$$
Meanwhile, integrating by parts in the first  integral by $t$ at the right-hand side of latter equality in (29) (the second one can be treated analogously) and using again (25), we find

$$\int_N^M  \cos (xt - \pi(2\nu+1)/4) \sqrt t \int_{\mu- i\infty}^{\mu+i\infty} g^*(s) \frac{\Gamma((1-s)/2)}{ \Gamma( (1- 2\nu + s)/2)} $$

$$\times Q^{-\nu}_{(s-1)/2} \left( {t^2+a^2 \over t^2-a^2}\right) 2^{-s}  \left(t^2- a^2\right)^{(s-1)/2} ds dt$$

$$=   O\left( N^{\mu+ {\rm Re}\nu- 1/2}\right) -  \int_N^M  \sin  (xt - \pi(2\nu+1)/4)  t^{3/2}  \int_{\mu- i\infty}^{\mu+i\infty} g^*(s) \frac{\Gamma((1-s)/2) (s-1) }{ \Gamma( (1- 2\nu + s)/2)} $$

$$\times Q^{-\nu}_{(s-1)/2} \left( {t^2+a^2 \over t^2-a^2}\right) 2^{-s}  \left(t^2- a^2\right)^{(s-3)/2} ds dt$$

$$-  \int_N^M  \sin  (xt - \pi(2\nu+1)/4)  t^{1/2}  \int_{\mu- i\infty}^{\mu+i\infty} g^*(s) \frac{\Gamma((1-s)/2) }{ \Gamma( (1- 2\nu + s)/2)} $$

$$\times {d\over dt} \left[ Q^{-\nu}_{(s-1)/2} \left( {t^2+a^2 \over t^2-a^2}\right)\right]  2^{-s}  \left(t^2- a^2\right)^{(s-1)/2} ds dt $$

$$=   O\left( N^{\mu+ {\rm Re}\nu- 1/2}\right) -    O\left( N^{\mu+ {\rm Re}\nu- 3/2}\right) $$

$$-  \int_N^M  \sin  (xt - \pi(2\nu+1)/4)  t^{1/2}  \int_{\mu- i\infty}^{\mu+i\infty} g^*(s) \frac{\Gamma((1-s)/2) }{ \Gamma( (1- 2\nu + s)/2)} $$

$$\times {d\over dt} \left[ Q^{-\nu}_{(s-1)/2} \left( {t^2+a^2 \over t^2-a^2}\right)\right]  2^{-s}  \left(t^2- a^2\right)^{(s-1)/2} ds dt, $$
where the differentiation under the integral sign is motivated by the absolute and uniform convergence and will be clearly seen form the estimate of the derivative of the second kind associated Legendre function, being deduced below. 
Indeed, recalling (23), we deduce

$${d\over dt}  \left[Q^{-\nu}_{ (s-1)/2} \left( t\right)\right]   =      \frac{ 2^{-(s+1)/2}   e^{-i\nu\pi} \nu \sqrt\pi\  \Gamma((s+1 - 2\nu)/2)}{ \Gamma((3+ s  -2 \nu)/4) \Gamma((1+s+2\nu)/4)}\  t  (t^2-1)^{-\nu /2}\left[ -  {\nu \over t^2-1} \right. $$

$$\times \int_0^1u^{- (1-s +2\nu)/4} (1-u)^{ (s-3+2\nu)/4} (t^2-u)^{ (2\nu-  1-s)/4} du$$

$$\left. +  { 2\nu-  1-s\over 2}  \int_0^1u^{- (1-s +2\nu)/4} (1-u)^{ (s-3+2\nu)/4} (t^2-u)^{ (2\nu-  5-s)/4}du\right],\ t > 1,$$
and, consequently,  in the same manner as in (25), we obtain 

$$\left| {d\over dt}  \left[Q^{-\nu}_{ (s-1)/2} \left( t\right)\right] \right| = O\left( |s|^{1- {\rm Re} \nu} t   (t^2-1)^{-\nu /2- 1}\right),\ t >1$$
under the same condition  $ {\rm Re} s > 2  {\rm Re} \nu -1$.   Hence, we easily get

$$\left|  {d\over dt} \left[ Q^{-\nu}_{(s-1)/2} \left( {t^2+a^2 \over t^2-a^2}\right)\right] \right| =
O\left( |s|^{1- {\rm Re} \nu} t   (t^2-a^2)^{\nu /2- 1}\right),\ t >a$$
and, finally, 

$$\left| \int_N^M  \sin  (xt - \pi(2\nu+1)/4)  t^{1/2}  \int_{\mu- i\infty}^{\mu+i\infty} g^*(s) \frac{\Gamma((1-s)/2) }{ \Gamma( (1- 2\nu + s)/2)} \right.$$

$$\left. \times {d\over dt} \left[ Q^{-\nu}_{(s-1)/2} \left( {t^2+a^2 \over t^2-a^2}\right)\right]  2^{-s}  \left(t^2- a^2\right)^{(s-1)/2} ds dt\right| \le C  \big\vert\big\vert g \big\vert\big\vert_{{\mathcal M}_{0,1}^{-1}(L_c)} \int_N^M  t^{  {\rm Re} \nu +\mu 
-3/2} dt$$

$$=   O\left( N^{\mu+ {\rm Re}\nu- 1/2}\right),\  N \to \infty,  \ \mu+ {\rm Re}\nu < 1/2. $$
Thus, combining with previous estimates, we see that $G_N$ is a Cauchy sequence and $x G_N(x)$ has the same pointwise limit (28) for almost all positive x.   Theorem is proved. 

\end{proof}

\bigskip
\centerline{{\bf Acknowledgments}}
\bigskip
The work was partially supported by CMUP [UID/MAT/00144/2013], which is funded by
FCT(Portugal) with national (MEC) and European structural funds through the programs FEDER,
under the partnership agreement PT2020.  The author thanks Mark Craddock for pointing out the Weber integral equation to his attention.  

\bigskip

\bibliographystyle{amsplain}

\begin{thebibliography}{}

\bibitem{erd}
 A. Erd\'elyi, W. Magnus, F. Oberhettinger and F.G. Tricomi,  {\it Higher Transcendental Functions}, Vols. I and  II, McGraw-Hill, New York, London and Toronto (1953).
 
 \bibitem {mar}  O.I. Marichev, {\it  Handbook of Integral Transforms of Higher Transcendental Functions: Theory and Algorithmic Tables},  Ellis Horwood,  Chichester  ( 1983).

\bibitem{prud}  A.P. Prudnikov, Yu. A. Brychkov and O. I.
Marichev,   {\it Integrals and Series: Vol. 1:  Elementary  Functions},
Gordon and Breach, New York  (1986); {\it Integrals and Series: Vol. 2: Special Functions},
Gordon and Breach, New York  (1986); {\it Vol. 3: More Special
Functions}, Gordon and Breach, New York  (1990).

\bibitem {titweb}  E.C. Titchmarsh,   Weber's integral theorem,  {\it Proc. Lond. Math. Soc.}, {\bf 22}  (2),  (1924),  15-28.

\bibitem {tit}  E.C. Titchmarsh, {\it  An Introduction to the Theory of Fourier Integrals}, Clarendon Press, Oxford ( 1937).

\bibitem {vira}  N. Virchenko, I.  Fedotova.  Generalized associated Legendre functions and their applications. (With a foreword by S. Yakubovich). World Scientific Publishing Co., Inc., River Edge, NJ,  (2001).

\bibitem {class} Vu Kim Tuan, O.I. Marichev and S. Yakubovich,  Composition structure of integral
transformations,  {\it J. Soviet Math.}, {\bf 33} (1986), 166-169.

\bibitem{yal}   S.  Yakubovich and  Yu.  Luchko, {\it The Hypergeometric Approach to Integral Transforms and Convolutions}.  Mathematics and its Applications, 287. Kluwer Academic Publishers Group, Dordrecht (1994).



\end{thebibliography}
{}

\end{document}